\documentclass{amsart}
\usepackage{color,amsmath,amssymb} 

\newcommand{\Cst}{\mathbb{C}^*}

\newcommand{\liea}{{\frak g}}
\newcommand{\csa}{{\frak h}}

\newcommand{\loopgr}{{\mathcal L}G }

\newcommand{\oloopgr}{{\mathcal L}_{Z}(\bar{G})}

\newcommand{\W}{W}
\newcommand{\Orb}{\mathcal{O}}
\newcommand{\C}{\mathbb{C}}
\newcommand{\Z}{\mathbb{Z}}
\newcommand{\N}{\mathbb{N}}

\newcommand{\Proj}{\mathbb{P}}

\newcommand{\g}{{\mathfrak{g}}}    
\newcommand{\h}{{\mathfrak{h}}}  
\newcommand{\n}{{\mathfrak{n}}}    

\renewcommand{\L}{{\mathcal{L}}}

\newcommand{\ad}{\textnormal{ad}}
\newcommand{\pol}{\textnormal{pol}}
\newcommand{\rk}{\textnormal{rk}}

\newtheorem{Theorem}{Theorem}[section]
\newtheorem{Lemma}[Theorem]{Lemma}
\newtheorem{Cor}[Theorem]{Corollary}
\newtheorem{Prop}[Theorem]{Proposition}
{\theoremstyle{definition}\newtheorem{Definition}[Theorem]{Definition}}
{\theoremstyle{remark}\newtheorem{Remark}[Theorem]{Remark}}
{\theoremstyle{remark}}

\begin{document}


\title[]{Conjugacy classes in Kac-Moody groups and Principal $G$-bundles over
  elliptic curves}


\dedicatory{\mbox{\normalfont\mdseries\small \today}}

\author{Stephan Mohrdieck and Robert Wendt}

\address{Stephan Mohrdieck\\Universit\"at Hamburg\\Mathematik\\
    Bundesstra\ss e 55\\20146 Hamburg\\Germany}
\email{stephan.mohrdieck@math.uni-hamburg.de}
\address{Robert Wendt\\University of Toronto\\Department of Mathematics\\
    100 St.George Street\\Toronto\\Ontario M5S 3G3\\Canada}
\email{rwendt@math.toronto.edu}



\setcounter{page}{1}

\begin{abstract}
For a simple complex Lie group $G$
the connected components of the moduli space of $G$-bundles
over an elliptic curve are weighted projective spaces.
In this note we will provide a new proof of this result using
the invariant theory of Kac-Moody groups, in particular the 
action of the (twisted) Coxeter element on the root system of 
$G$.
\end{abstract}

\maketitle

\section{Introduction}
Let $G$ be a simple algebraic group over $\C$. It is known (\cite{L, FM2}) 
 that the connected components of the moduli space of semistable
 $G$-bundles over an elliptic curve are isomorphic to  weighted
 projective spaces. 
 In this note, we show how this result can be obtained from the geometry
 of the holomorphic Kac-Moody group $\widetilde G$ associated to $G$. 
The main step is a detailed description of the action of the (twisted) 
Coxeter element of the Weyl group of $\widetilde G$ on the root system of $G$. In particular,
we show that the holomorphic principal bundle associated to the Coxeter element is minimally 
unstable in the sense of \cite{FM2}.

Let us briefly sketch the main idea assuming that $G$ is simply connected. 
Let $\L(G)$ denote the holomorphic loop group of $G$, 
i.e. the group of holomorphic maps from $\C^*\to G$. Fix some
$q\in\C^*$ with $|q|<1$. 
It has been
observed by E.Looijenga (\cite{EF, BG})  that there is a bijection between the set of 
$q$-twisted conjugacy classes in the group $\L(G)$ and
the set of isomorphism classes of holomorphic $G$-bundles on $E_q=\C^*/q^\Z$.
This observation suggests to use geometric
invariant theory for $\L(G)$ to describe the moduli space of holomorphic
principal $G$ bundles on $E_q$. However, in order to obtain a good 
invariant theory we have to pass to a central extension of the group $\L(G)$,
i.e. to the holomorphic Kac-Moody group $\widetilde G$ corresponding to $G$.
Br\"uchert \cite{Br} has
constructed an  analogue of the Steinberg cross section in $\widetilde G$. 
This cross section
carries a natural $\C^*$-action so that the quotient is isomorphic to a
weighted projective space. Since the Kac-Moody group $\widetilde G$ is a
central extension of the loop group 
$\L(G)$, we can associate to each element of the cross
section a holomorphic $G$ bundle. It turns out that bundles in the same
$\C^*$-orbit are isomorphic and that the space of $\C^*$-orbits is isomorphic
to the moduli space of semistable $G$ bundles over the elliptic curve.
The main step here is to show that the bundles coming from the cross section
are semistable outside of the section's origin and hence have an image in the 
moduli space. This is carried out as follows:
The origin of the section corresponds to the (twisted) Coxeter element.
A careful investigation of its action on the set of all roots of $G$
shows that the corresponding principal bundle is unstable and has minimal
possible automorphism group dimension among all unstable bundles 
(Proposition \ref{stablesec}). Since this bundle is degeneration of all
the bundles corresponding to other elements of the cross section,
a result of Helmke and Slodowy \cite{HS} implies
that the latter bundles have smaller automorphism group dimension and
hence are semistable. 

The ideas here are mostly due to Peter Slodowy. The Steinberg cross
section in $\widetilde G$ plays a role in a generalisation of a theorem 
of Brieskorn
which relates simple singularities of type A$_n$, D$_n$ and E$_n$ 
and the corresponding simple 
Lie groups to the case of elliptic singularities (see \cite{HS1}).

\medskip

The paper is set up as follows: In section \ref{KM} we recall some results 
from the Theory of affine Kac-Moody Lie algebras and groups. 	n particular, 
we describe the Steinberg cross section and mention how some 
standard results on Kac-Moody groups generalise to the case when the 
underlying finite-dimensional Lie group is not simply connected. In section 
\ref{holbundles} we make the connection to holomorphic bundles.  Section \ref{sec:proof} 
we give a proof of Proposition \ref{stablesec} while the Appendix gives the details 
of some explicit calculations.

\medskip

{\it Acknowledgement:} The first author expresses his gratitude to the DFG for financial support.


\section{Holomorphic Kac Moody Groups}\label{KM}
\subsection{The group}\label{group}
Let $G$ be a simple and simply connected algebraic group over $\C$ and denote by  
$\L(G)$ group of holomorphic maps from $\C^*\to G$ endowed with point-wise multiplication.
This is an infinite-dimensional Lie group, called the holomorphic loop group corresponding to $G$.
The group $\L(G)$
possesses a universal central extension which sits in an exact sequence
$$
  1\to \C^* \stackrel{\iota}{\longrightarrow} \widehat \L(G) 
  \stackrel{\pi}{\longrightarrow} \L(G) \to 1\,.
$$
Topologically, the group $\widehat \L(G)$ is a non-trivial $\C^*$-bundle over
$\L(G)$. Taking the $l$-th power of the corresponding transition functions we obtain the
central extension of $\L(G)$ of level $l$.Up to isomorphism this yields all central 
extensions of $\L(G)$ (\cite{PS}).

The natural multiplication action of the multiplicative group
$\C^*$ on the loop group $\L(G)$ lifts to a $\C^*$-action on 
the central extension $\widehat \L(G)$, and we define
the holomorphic Kac-Moody group corresponding to $G$ to be the semi-direct product
$\widetilde G=\widehat \L(G)\rtimes\C^*$.
The Lie algebras of $G$, $\L(G)$ and $\widetilde G$ are denoted by $\g$,
$\L(\g)$ and $\widetilde \g$, respectively.
 

\subsection{Roots and reflections}
If the finite-dimensional Lie algebra $\g$ is simple of rank $r$, the subalgebra
$\widetilde\g_{\pol}=\g\otimes\C[z,z^{-1}]\oplus\C C\oplus\C D
\subset\widetilde \g$ of polynomial loops is an untwisted 
affine 
Lie algebra in the sense of  \cite{K}, and $\widetilde \g$ 
can be viewed as 
a certain completion of it (see \cite{GW}). Here, $C$ denotes a generator of the centre of 
$\widetilde\g_{\pol}$, and $D$ is the infinitesimal generator of the $\C^*$-action on the centrally 
extended group $\widehat\L(G)$.

Let us fix once and for all a maximal torus $T\subset G$, and denote the corresponding Lie algebra 
$Lie(T)=\mathfrak h\subset \g$.
Is known that
the Lie algebra $\widetilde \g_{\pol}$ has a root space decomposition with respect to the Cartan subalgebra 
$\h\oplus\C C\oplus \C D\subset\widetilde \g_{\pol}$.
Denote the set of roots by $\widetilde\Delta$ and fix a set
$\widetilde{\Pi}=\{\widetilde{\alpha}_0,\ldots,\widetilde{\alpha}_r\}$  of simple roots.
We get a linear combination $C=\sum_{i=0}^ra_i^\vee\widetilde{\alpha}_i^\vee$ 
where the 
$\widetilde{\alpha}_i^\vee $ are the simple co-roots. The $a_i^\vee$
appearing in the expression above are called the dual Kac labels.

Denote the by $\widetilde W$ the Weyl group of the root system $\widetilde\Delta$ and 
let $r_i$ be the
simple reflection corresponding to $\widetilde\alpha_i\in\widetilde\Pi$.
The Weyl group $\widetilde W$ is known to have the following two descriptions:
\begin{eqnarray*}
\widetilde W&\cong &W\ltimes Q^\vee\,,\\
\widetilde W&\cong &
  N_{\L(G)\rtimes\C^*}(T\times\C^*)\,/\,(T\times\C^*)\,.
\end{eqnarray*}
Here, $W$ is the Weyl group of the group $G$ and $Q^\vee$ is the co-root lattice of $\g$.
Furthermore, $N_{\L(G)\rtimes\C^*}(T\times\C^*)$ denotes the normaliser
of the torus $T\times\C^*$ in $\L(G)\rtimes\C^*$, (see \cite{PS}).

Finally, the product of all simple reflections $cox=\prod_{i=0}^rr_i$ is called a Coxeter element of
$\widetilde W$. Obviously this definition depends on the choice of a basis of the
root system $\Delta$. But it is known (see e.g. \cite{Hu}) that all  Coxeter elements of $\widetilde W$ are
conjugate in $\widetilde W$ unless $\g$ is
of type $A_n$.


\subsection{Representations and characters}
Denote by $\widetilde{P}$ the weight lattice of $\widetilde\g_{\pol}$
and by $\widetilde{P}^+$ its cone of dominant weights which is generated by the 
fundamental dominant weights $\lambda_0,...,\lambda_r$ and $\delta$ (throughout this section, 
we follow the name conventions of \cite{K}).
For each $\lambda\in \widetilde{P}^+$ there is an irreducible highest weight
module $V_{\lambda}$ of $\widetilde\g_{\pol}$.
It is known that each of these representations 
extends to a representation of $\widetilde\g$ on the analytic
completion $V_{\lambda}^{an}$
of $V_{\lambda}$ which, in turn, lifts to a representation
of the holomorphic Kac-Moody group $\widetilde G$ (see \cite{GW}). 

The vector space $V_{\lambda}$ admits a positive definite
Hermitian form $(.,.)$ which is contravariant with respect to the anti-linear
Cartan involution on $\widetilde\g_{\pol}$. We 
denote by $V^{ss}_{\lambda}$ the $L^2$-completion of $V_{\lambda}$
with respect to the norm defined by this Hermitian form.

For a fixed $q\in\C^*$, let the ``$q$-level set'' of $\widetilde G$ be the subset
$\widetilde G_q= \widehat \L(G)\times\{q\}\subset\widetilde G$.
Obviously, the $q$-level sets $\widetilde G_q$ are invariant under conjugation in
$\widetilde G$.

The following
Theorem is known (see \cite{GW}, \cite{EFK}, \cite{B})
\begin{Theorem}\label{trace} Fix $q\in\C^*$ such that $|q|<1$. Then we have
 \renewcommand{\theenumi}{\roman{enumi}}
\renewcommand{\labelenumi}{\textnormal{(\theenumi)}}
\begin{enumerate}
\item
  For any $(g,q)\in\widetilde G_q$, 
  the operator
  $(g,q):V^{an}_{\lambda}\to V^{an}_{\lambda}$ uniquely   
  extends to a trace class operator on $V^{ss}_{\lambda}$.
\item
  The function 
  $\chi_{\lambda}:\widetilde G_q\to\C$, defined by
  $(g,q)\mapsto tr|_{V^{ss}_{\lambda}}(g,q)$
  is holomorphic and conjugacy invariant.
\item
  The functions 
  $\chi_{\lambda_0},\ldots,
  \chi_{\lambda_r}$ generate the ring of holomorphic conjugacy
  invariant functions on $\widetilde G_q$.
\end{enumerate}
\end{Theorem}
Note that the centre of $\widetilde{\g}$ acts on $V_{\lambda_i}$ by the scalar $a_i^\vee$.
Hence, for any $q\in\C^*$ such that $|q|<1$, we can define 
a conjugacy invariant
map $\chi_q:\widetilde G_q \to \C^{r+1}$ via 
$(g,q)\mapsto
(\chi_{\lambda_0}(g,q),\ldots\chi_{\lambda_r}(g,q))$.


\subsection{The cross section}\label{cross}
In this section we shall indicate, following
\cite{Br}, how to define a section to the map $\chi_q$.

Associated to each real root $\widetilde\alpha \in\widetilde\Delta$, (i.e. a root whose $\widetilde W$-orbit
$\widetilde W \widetilde\alpha $ contains a simple root,)
there exists a one parameter subgroup
$x_{\widetilde\alpha}: \C\to \widehat \L(G)$ (see \cite{GW}).
For $0\leq i\leq r$ and $c\in\C$, we set $x_i(c)=x_{\alpha_i}(c)$, 
$y_i(c)=x_{-\alpha_i}(c)$, and $n_i= x_i(1)y_i(-1) x_i(1)$.
(Thus the element $n_i\in\widehat\L(G)$ 
is a representative of the simple reflection $r_i$ in the Weyl group $\widetilde W$.)

For each $q\in\C^*$  we can define a map 
$\omega_q:\C^{r+1}\to \widetilde G_q$ via
$$
  \omega_q:
  (c_0,\ldots,c_r) \mapsto \left(\prod_{i=0}^r x_i(c_i)n_i,q\right)\,.
$$
The map $\omega_q$ as well as its image 
$\mathcal C_q=\omega_q(\C^{r+1})$ are called a 
Steinberg cross-section in
$\widetilde G_q$. This name is justified by the analogy of the
definition of the map $\omega_q$ to the usual Steinberg cross section
for finite dimensional algebraic groups. The following Theorem is due
to Br\"uchert (\cite{Br}, Theorem 7):
\begin{Theorem}\label{Bru1}
  For $|q|$ small enough, the map
  $\chi_q\circ\omega_q:\C^{r+1}\to\C^{r+1}$ is an isomorphism of
  algebraic varieties.
\end{Theorem}

An essential tool in the proof of Theorem \ref{Bru1} is the existence
of a $\C^*$-action on the cross section  $\mathcal C_q$. Let
$\iota:\C^*\to\widehat \L(G)$ denote the identification of the
centre of $\widehat \L(G)$ with $\C^*$. Then we have (\cite{Br},
Proposition 15):

\begin{Prop}\label{Bru2}
There exists a one parameter subgroup $\mu:\Cst\to \widehat \L(G)$
such that for any $(c,q)\in \mathcal C_q$, we have 
$(\mu(z)\,c\,\iota(z)\mu(z)^{-1},q)\in\mathcal C_q$, for all $z\in\Cst $.
\end{Prop}

Using Proposition \ref{Bru2} we can define a $\C^*$-action on
$\mathcal C_q$ via 
$$
  z:(c,q)\mapsto (\mu(z)\,c\,\iota(z)\mu(z)^{-1},q)\,.
$$
The trace function $\chi_q$ behaves well with respect to this
$\C^*$-action. In fact, we have (\cite{Br}, Theorem 6):

\begin{Prop}  
  There exists some $k\in\N$ such that 
  $\chi_i(z.(c,q))=z^{ka_i^\vee}\chi_i(c,q)$ for all 
  $(c,q)\in\mathcal C_q$.
\end{Prop}

\begin{Cor}\label{WP1}
  If we let $\C^*$ act on $\C^{r+1}$ with weights $ka_i^\vee$,
  then the restriction of  $\chi_q$ to $\mathcal C_q$ induces a
  $\C^*$-equivariant isomorphism $\mathcal C_q\to\C^{r+1}$. In
  particular, we have
  $$
    \left(\mathcal C_q - \{\omega(0)\}\right) /\C^* 
     \cong \Proj(a_0^\vee,\ldots a_r^\vee)\,,
  $$
  where $\Proj(a_0^\vee,\ldots a_r^\vee)$ denotes the weighted projective space
  with weights $a_0^\vee,\ldots a_r^\vee$.
\end{Cor}


\subsection{The non-simply connected case}\label{nonsimply}
Hewe we give a brief account on how the constructions of the previous paragraphs
generalise to the non simply connected case. A more thorough treatment
can be found in \cite{M}.

Let $G$ be a simple algebraic group with universal cover $\bar{G}$ with 
maximal torus $T$ resp. $\bar{T}$ and co-character lattices
$\check{\chi}(T)$ and $\check{Q}=\check{\chi}(\bar{T})$.
Denote by $Z\cong \pi_1(G)$ the kernel of the covering map $\bar{G}\to G$.

Then the component group of $\L(G)$ is given by $\pi_0({\mathcal L}(G))\cong Z$.

For the representation theory it turns out to be more appropriate to work with 
the ''group of open loops'' instead. We set:
$$
\oloopgr =
\{\phi: \C\to \bar{G}\text{ such that } \phi(t)\phi(t+1)^{-1}\in Z,\,\text{ and }
\phi\, \text{ holomorphic}\}\,.
$$
Identifying $\check{\chi}(T)$ and $\check{Q}$ as lattices in $\csa $, the 
exponential map provides
a group homomorphism:
\begin{eqnarray*}
\gamma: \check{\chi}(T)&\to&\oloopgr\\ 
\check{\beta}&\mapsto&\gamma_{\check{\beta}}:t\mapsto
\gamma_{\check{\beta}}(t):=e^{2\pi i \check{\beta} t}.
\end{eqnarray*}
With these notations the following result holds:\newpage
\begin{Lemma}\quad
\renewcommand{\theenumi}{\roman{enumi}}
\renewcommand{\labelenumi}{\textnormal{(\theenumi)}}
\begin{enumerate}
\item
$\gamma(\check{Q})\subset \oloopgr$,
\item
$\oloopgr \cong ({\mathcal L}\bar{G}\rtimes\check{\chi}(T))/\check{Q}$,
\item
$\loopgr\cong\oloopgr/Z$.
\end{enumerate}
\end{Lemma}

The group $Z\subset G$ can be canonically identified with a subgroup of the group of diagram 
automorphisms on the Dynkin diagram of the root system $\widetilde \Delta$ (see e.g. \cite{B, T}).
Since we are interested in a cyclic component group
let us fix an automorphism $\sigma\in Z$ and denote 
by $\Sigma$ the subgroup 
generated by the element $\sigma$. Choose a representative  
$\check{\lambda}_{\sigma}\in \check{\chi}(T)$ of $\sigma$ and write $\bar{\Sigma}\cong\Z$
for the subgroup (of $\check{\chi}(T) $) generated by $\check{\lambda}_{\sigma} $.
One finds $(\check{\lambda}_{\sigma})^{\mbox{\tiny ord}\sigma}\in \check{Q}$
Consider the group 
$$
{\mathcal L}_{\Sigma}\bar{G}:=
({\mathcal L}\bar{G}\rtimes\bar{\Sigma})/\bar{\Sigma}^{\mbox{\tiny ord}\sigma}.
$$
In \cite{M}, Theorem 3.2, see also \cite{T}, Section 3.3, there is a 
classification of 
all central extensions of ${\mathcal L}_{\Sigma}\bar{G}$:
\begin{Theorem}
\quad
\renewcommand{\theenumi}{\roman{enumi}}
\renewcommand{\labelenumi}{\textnormal{(\theenumi)}}
\begin{enumerate}
\item
There is a natural number $k_f\in \N$ such that for all $l\in k_f\Z$  
a uniquely determined 
central extension $\widehat{\mathcal L}_{\Sigma}\bar{G}^l$ 
of ${\mathcal L}_{\Sigma}\bar{G}$
of level $l$ exists.
\item 
Every central extension of  
${\mathcal L}_{\Sigma}\bar{G}$ is isomorphic to $\widehat{\mathcal L}_{\Sigma}\bar{G}^l$ for a certain $l$.\\
\item
The translation action of $\C$ on ${\mathcal L}_{\Sigma}\bar{G}$ 
lifts to these central extensions and it
factors through $m\Z$ for a certain $m\in\N$.
\end{enumerate}
\end{Theorem}
Writing $\widetilde{\Cst}$ for $\C/m\Z$  this allows us to define the 
non-connected Kac-Moody group:
\begin{Definition}
The group $\tilde{G}:=\widehat{{\mathcal L}}_{\Sigma}\bar{G}^{k_f}\rtimes\widetilde{\Cst}$ 
is called the affine Kac-Moody group associated to $G$ and $\Sigma$.
\end{Definition} 
If $\lambda$ is a $\Sigma$-invariant dominant weight then its level $k$ (the multiple
by which the centre acts on the representation space) is divisible by $k_f$.
Furthermore the action of $\widehat{{\mathcal L}}_{\Sigma}\bar{G}_0$ extends to a 
$\widehat{{\mathcal L}}_{\Sigma}\bar{G}$-module with Theorem 2.1 (i) and (ii) still being valid for 
$(g,q)\in \widehat{{\mathcal L}}_{\Sigma}\bar{G}_0$ with 
$|q|<1$. 
For $s+1$ being the number of $\Sigma $-orbits on $\widetilde{\Pi}$
there are $s+1$ dominant weights $\Lambda_0,...,\Lambda_s$ with
$\widetilde{P}^{+,\Sigma }:=\{\lambda\in \widetilde{P}^{+},\sigma (\lambda)=\lambda \}=
\Z^{\geq 0}<\Lambda_0,...,\Lambda_s>\times\Z\delta $. 
Write $\tilde{G}_{\sigma}$ for the connected component of 
$\tilde{G}$  corresponding to $\sigma $. Then the 
following analogue of Theorem 2.1(iii) is valid:
\begin{Prop}
The functions $\chi_{\Lambda_0},...,\chi_{\Lambda_s}$ generate the ring 
of holomorphic conjugacy invariant 
functions on $\tilde{G}_{\sigma, q}$, 
for $q\in D^{\ast}$.
\end{Prop}
Again, this allows us to define the quotient map $\chi_q:
\tilde{G}_{\sigma, q}\to 
\C^{s+1}$ via $(g,q)\mapsto (\chi_{\Lambda_0}(g,q),...,\chi_{\Lambda_s}(g,q))$.

It is also possible to construct a cross section to this quotient map in this 
situation:
Choose representatives $\{\tilde{\alpha }_0,...,\tilde{\alpha }_s\}
\subset\widetilde{\Pi}$ of the $\Sigma$-orbits
on $\widetilde{\Pi}$ and a representative $n_{\sigma}\in 
N_{\bar{G}}(\hat{T})$ of $\sigma$. 
Here, $N_{\bar{G}}(\hat{T})$ is the normaliser of $\hat{T}$ in $\bar{G}$.
Using the notation of section 2.5, define a map:
\begin{eqnarray*}
\omega_{\sigma, q}:\C^{s+1}&\to&\tilde{G}_{\sigma, q}\\
\omega_{\sigma, q}(c_0,...,c_s)&:=&\left(\left(\prod_0^s x_i(c_i)n_i \right)
n_{\sigma},q\right).
\end{eqnarray*}
We set $\mathcal C_{\sigma, q}:=\omega_{\sigma, q}(\C^{s+1})$.
Note that $\omega_{\sigma, q}(0,...,0)$ is a representative of the twisted
Coxeter element 
$\mbox{cox}^{\sigma}=s_0...s_s\sigma\in\widetilde{W}\rtimes\pi_1(G)$.
(For the definition of twisted Coxeter elements, see \cite{Sp,M}.)

It is shown in \cite{M}, Theorem 4.2, Lemma 4.1 and Lemma 2.9, that that the results 
Theorem 2.2, Propositions
2.3 an 2.4 and Corollary 2.5 also carry over to this setting, if we replace the 
dual Kac labels $a_i^\vee$ 
by the ones of the orbit Lie algebra in the sense of Fuchs etal. \cite{FSS}.
In particular, we have:
\begin{Prop}
\label{section}
 For $|q|$ small enough, the map
  $\chi_{\sigma, q}\circ\omega_{\sigma, q}:\C^{s+1}\to\C^{s+1}$ is an isomorphism of
  algebraic varieties.
\end{Prop}



\section{Holomorphic principal G-bundles } \label{holbundles}
\subsection{Gluing maps for $G$-bundles}\label{glue}
In this section, we relate the conjugacy classes in 
the $q$-level set $\widetilde G_q \subset \widetilde G$ 
to principal $G$-bundles over the elliptic curve $E_q=\C^*/q^\Z$.
Up to
$C^\infty$-isomorphism,
every principal $G$-bundle over $E_q$ is determined by its topological class,
which is an element in $\pi_1(G)\cong Z$. One can classify holomorphic
principal $G$-bundles of a fixed topological class $\sigma\in Z$ as follows. 
We have $\pi_0(\L(G))\cong Z$. 
Hence we can consider the connected component $\L(G)_{\sigma}$ of the loop group 
$\L(G)$ which corresponds to the element $\sigma\in Z$. 
The group $\L(G)\rtimes\C^*$ acts on the set
$\L(G)_{\sigma}\times\{q\}\subset \L(G)\rtimes\C^*$ by conjugation.
A fundamental observation due to E. Looijenga gives a one-to-one
correspondence between the set of equivalence classes of
holomorphic $G$-bundles on $E_q$ of
topological type $\sigma$ and the set of $\L(G)\times\{1\}$-orbits in
$\L(G)_{\sigma}\times\{q\}$. This correspondence comes about as follows. For any
element 
$(\gamma,q)\in \L(G)_{\sigma}\times\{q\}$ consider the 
$G$-bundle $\xi_\gamma$ over $E_q$
which is defined as  the quotient $(\C^*\times G)/\Z$, where $\Z$ acts via
$1:(z,h)\mapsto (qz,\gamma(z)h)$. Obviously, this construction defines
a holomorphic $G$-bundle $\xi_\gamma$
of topological type $\sigma$ on the elliptic curve
$E_q$, and the following result due to E. Looijenga 
is not hard to prove (see e.g. \cite{EF, BG}).
\begin{Theorem}\label{bundles}\quad
\renewcommand{\theenumi}{\roman{enumi}}
\renewcommand{\labelenumi}{\textnormal{(\theenumi)}}
\begin{enumerate}
\item
  Two elements $(\gamma_1,q),\,(\gamma_2,q)\in \L(G)_{\sigma}\times\{q\}$ 
  are conjugate under
  $\L(G)\times\{1\}$ if and only if the corresponding holomorphic
  $G$-bundles $\xi_{\gamma_1}$ and
  $\xi_{\gamma_2}$ are isomorphic.  
\item
  For any holomorphic 
  $G$-bundle $\xi\to E_q$ of topological type $\sigma$, there
  exists an element $\gamma\in \L(G)_{\sigma}$ such that $\xi\cong \xi_\gamma$.
\end{enumerate}
\end{Theorem}


\subsection{Vector bundles and stability}
Here, we recall the notion of stability of principal bundles over $E_q$ and the reduction
of unstable bundles to a Levi subgroup which corresponds to the 
Harder-Narashiman filtration in the vector bundle case.
The slope $\mu(V)$ of a holomorphic vector bundle $V\to E_q$ over the
elliptic curve $E_q$ (or over any smooth curve) is defined by 
$\mu(V)=\deg(V)/\rk(V)$. The vector bundle $V$ is said to be stable
(resp. semistable) if 
$$
  \mu(U) < \mu(V)\qquad \left(\text{resp.}\quad \mu(U) \leq \mu(V)\right)
$$
for any nontrivial holomorphic subbundle $U$ of $V$.

\begin{Definition} 
  A principal bundle $\xi\to E_q$ is called semistable if the
  associated vector bundle $\ad(\xi)=\xi\times^G\g$ 
  is a semistable vector bundle. The
  principal bundle $\xi\to E_q$ is called unstable if it is not
  semistable.
\end{Definition}

For any holomorphic vector bundle $V$ on $E_q$ there exists a
unique filtration, the so called Harder Narashiman filtration, 
$0=V_0\subset V_1\subset\ldots\subset V_n=V$ such
that the quotients $V_i/V_{i-1}$ are semistable and 
$\mu(V_i/V_{i-1})>\mu(V_{i+1}/V_{i})$. Obviously, the bundle $V$ is semistable
exactly if $n=1$. Since $E_q$ is an elliptic curve, this filtration
splits so that we have
$$
  V\cong\bigoplus V_i/V_{i-1}\,.
$$
The existence of the filtration 
$0=V_0\subset V_1\subset\ldots\subset V_n=V$
corresponds to a reduction of the structure group of
the vector bundle $V$ from $GL(V)$ to a parabolic subgroup 
$P\subset GL(V)$, and the splitting of the filtration 
corresponds to a further
reduction of the structure group to a Levi subgroup $L$ of $P$.
    
This construction generalises to holomorphic $G$-bundles as follows
(see e.g. \cite{FM1}, Theorem 2.7., \cite{HS}, Theorem 1.3.1)
\begin{Theorem}\label{red}
  Let $\xi$ be an unstable holomorphic principal $G$-bundle over the
  elliptic curve $E_q$. Then there exists a maximal parabolic subgroup
  $P$ of $G$ and a Levi subgroup $L$ of $P$ such that the bundle
  $\xi$ reduces to
  a semistable $L$-bundle $\xi_L$ such that for any nonzero dominant
  character $\psi:P\to\C^*$, the line bundle associated to $\xi_L$ and
  $\psi$ has positive degree.
\end{Theorem}


\subsection{Bundles for elements in the Steinberg cross section}
Consider the $q$-level set  $\widetilde G_q\subset\widetilde G$ for some fixed $q$ with $|q|<1$.
The projection $\pi_q:\widetilde G_q\to \L(G)\times\{q\}$,
$(g,q)\mapsto(\pi(g),q)$ 
allows to associate to each element $(g,q)\in \widetilde G_q$ the
$G$-bundle on the elliptic curve $E_q$ with gluing map $\pi(g)$.
Fix an element $\sigma\in\pi_1(G)$. In this section we determine those $(g,q)\in\mathcal C_{\sigma,q}$ 
of the Steinberg 
section whose corresponding principal bundle $\xi_{\pi(g)}$ (which is of topological type $\sigma$) 
is semistable.
 
Assume $\pi_1(G)$ to be cyclic, possibly trivial, with generator $\sigma$ and let $s+1$ be 
the number of $\sigma$-orbits on the affine Dynkin diagram of $\widetilde \Delta$. 
 Then the following statement holds:
\begin{Prop}\label{prop:cox}
If $\xi_{\sigma}$ is the $G$-bundle of topological type $\sigma$
over $E_q$ corresponding to the twisted Coxeter element 
then $\mbox{dim}\,\mbox{Aut}(\xi_{\sigma})=s+1$.
\label{stablesec}
\end{Prop}
We give a proof of this Proposition in section \ref{sec:proof}.
\begin{Remark}
Using the terminology of \cite{FM2} we have shown that the twisted Coxeter element gives rise to 
a minimally unstable bundle.
\end{Remark}
As before, we denote the Steinberg cross section by $\mathcal C_{\sigma,q}=\omega_{\sigma,q}(\C^{s+1})$ and set 
$\mathcal C_{\sigma,q}^*=\mathcal C_{\sigma,q}\backslash\{\omega_{\sigma,q}(0)\}$. Note that 
Proposition \ref{prop:cox} implies the following corollary:
\begin{Cor}\quad
 \renewcommand{\theenumi}{\roman{enumi}}
\renewcommand{\labelenumi}{\textnormal{(\theenumi)}}
\begin{enumerate}
\item
The bundle $\xi_{\sigma}$ is unstable with minimal possible automorphism group dimension.
\item
All bundles corresponding to elements in $\mathcal C_{\sigma,q}^*$ are semistable.
\end{enumerate}
\end{Cor}
\begin{proof}
Using the construction from
section \ref{glue}, we associate to each element of the Steinberg cross section
$\mathcal C_{\sigma, q}$
a holomorphic $G$-bundle on $E_q$ of topological type $\sigma$. 

Consider the following family $\widetilde{\Xi}$ of $G$-bundles over 
the pointed complex space $(\mathcal C_{\sigma,q},\omega_{\sigma,q}(0))$:
$$
\widetilde{\Xi} = (\C^*\times \mathcal C_{\sigma,q}\times G)/\Z\,,
$$
with $1\in\Z$ acting by:
$$
  1:(z,(g,q),h)\mapsto (qz,(g,q),\pi(g)(z)h)\,.
$$
As usual, $\pi$ denotes the projection of the centrally extended loop group to the loop group itself.
The special fibre of this family is the bundle $\xi_{\sigma}$ from Proposition \ref{prop:cox}
and the construction shows that this 
bundle is adjacent to all the other bundles $\eta$ arising from elements in the Steinberg section 
$\mathcal C_{\sigma,q}^*$.
(A bundle $\xi$ is said to be adjacent to a bundle $\eta$
if there is a family of $G$-bundles over a pointed complex space $(\mathcal B,b)$ with
$\xi$ being isomorphic to the special fibre and each other fibre 
being isomorphic to $\eta$.)
Here we use that the section admits a $\Cst$-action with $\omega_{\sigma,q}(0)$ being contained
in every orbit closure.
By \cite{HS} Proposition 3.9 this implies $\mbox{dim}\mbox{Aut}(\eta)\leq s$. In view of 
\cite{HS}, Proposition 3.3 and Theorem 2.2, the each bundle $\eta\in \mathcal C_{\sigma,q}^*$
has to be (regular) semistable.
\end{proof}

\begin{Remark}
Using Atiyah's classification of vector bundles on an elliptic curve, one can show that if the group $G$ 
is simply connected then the $G$-bundle $\xi_g$ corresponding to an element $(g,q)\in\widetilde G_q$ 
is unstable exactly if $\chi_q(g,q)=0$. In particular, the bundles corresponding to
Weyl group elements of infinite order are unstable.
Presumably this statement also holds for the non-simply connected case.
\end{Remark}


\subsection{The moduli space of holomorphic G-bundles}
In this section we aim to describe the coarse 
moduli space $\mathcal M_{\sigma ,q}(G)$ of semistable
$G$-bundles over the elliptic curve $E_q$ of topological type $\sigma\in\pi_1(G) $
where we assume that $\sigma $ generates $\pi_1(G) $. 
We keep the notation of the previous section.
If $(g_1,q)$ and $(g_2,q)$ are two elements
of $\mathcal C_{\sigma ,q}^*$ which lie in the same $\C^*$-orbit on 
$\mathcal C_{\sigma ,q}^*$, then the bundles corresponding to $(g_1,q)$ and $(g_2,q)$ 
are isomorphic. 
In the non-simply connected case keep in mind
that the loop group $\loopgr$ is obtained from the group of open loops
$\oloopgr $ by further division by $Z=\pi_1(G)$.
However, a computation using the isomorphism given by Proposition \ref{section} shows
that the product $(cg,q) $ of an element $c\in Z$ and an element $(g,q)\in\mathcal C_{\sigma ,q}^*$
is conjugate to an element in the $\Cst$-orbit of $\xi$. 

Therefore we obtain a map 
$$
  \widetilde\psi:\mathcal C_{\sigma,q}^*\to \mathcal M_{\sigma,q}(G)\,,
$$
which factors through an injective map 
$$
   \psi:\mathcal C_{\sigma,q}^*/\C^* \to \mathcal M_{\sigma,q}(G)\,.
$$
It remains to show that $\psi$ is an algebraic isomorphism.
To this end, let us construct a holomorphic $G$-bundle $\Xi$ on
$E_q\times\mathcal C_{\sigma,q}^*$ as the quotient
$$
  \Xi = (\C^*\times \mathcal C_{\sigma,q}^* \times G)/\Z\,,
$$
where $\Z$ acts via
$$
  1:(z,(g,q),h)\mapsto (qz,(g,q),\pi(g)(z)h)\,.
$$
This shows that the map $\mathcal C_{\sigma,q}^*\to\mathcal M_{\sigma,q}(G)$ is
algebraic. Consequently, the induced map 
$\psi:\mathcal C_{\sigma,q}^*/\C^* \to\mathcal M_{\sigma,q}(G)$
is algebraic as well.

Now, $\psi$ is an injective morphism of irreducible
projective varieties of the same dimension. So it has to be an
isomorphism.
Hence we have proved
\begin{Theorem} 
  Let $\mathcal C_{\sigma,q}^*=\omega_{\sigma,q}(\C^{s+1}-\{0\})$ denote the cross
  section in $\widetilde G_{\sigma,q}$ and let $\mathcal M_{\sigma,q}(G)$ denote the moduli
  space of semistable $G$-bundles on the elliptic curve $E_q$. Then  
  $$
    \mathcal M_{\sigma,q}(G)\cong \mathcal C_{\sigma,q}^*/\C^*\,.
  $$
\end{Theorem}
In particular, using Corollary \ref{WP1}, respectively its analogue in the non-simply connected case
(see section \ref{nonsimply})
we get the following well
known result (see e.g. \cite{FM2}, \cite{L} and \cite{S} for the non-simply connected case).
\begin{Cor} 
Let $a_0^\vee,\ldots,a_s^\vee$ denote the dual Kac-labels of $G$ (respectively, the dual
Kac-labels of the orbit Lie algebra of $G$) for $\sigma\neq id$. Then
$$
  \mathcal M_{\sigma,q}(G)\cong\Proj(a_0^\vee,\ldots a_s^\vee)\,.
$$
\end{Cor}



\section{Proof of Proposition \ref{stablesec}}\label{sec:proof}
\begin{proof}
We have to show that the automorphism group $\mbox{Aut}(\xi_{\sigma})$ of the bundle corresponding to the twisted
Coxeter element has dimension $s+1$.
This is carried out by a case-by-case analysis. 
The corresponding calculations will be carried out over the adjoint bundle
$\mbox{ad}(\xi_{\sigma})=\xi_{\sigma}\times^{G}\liea $.
It is easy to see that $\mbox{ad}(\xi_{\sigma})$ can be described as the quotient
$(\Cst\times\liea )/\Z$ where $1\in\Z$ acts by $(z,X)\mapsto(zq,\mbox{Ad}(\gamma_{\sigma}(z))(X))$.
Here $\gamma_{\sigma}$ is the open loop corresponding to the twisted Coxeter element.

{\bf Step 1:} We first determine the Harder-Narasimhan filtration of the bundle $\xi_\sigma$. 
This is achieved by 
considering the action of $\gamma_{\sigma}$ on the Lie algebra. One easily sees that 
subbundles of $\mbox{ad}(\xi_{\sigma})$ correspond to subspaces of $\liea$ which are invariant 
under $Ad(\gamma_{\sigma}(z))$ for all $z\in\Cst $.
After having fixed a basis for the root system, the loop $\gamma_{\sigma}$ has the form:
$\gamma_{\sigma}=n_0...n_s\check{\lambda}(z)n_w$ for our representatives of the simple reflections $s_i$
(with exactly one reflection chosen for each $\sigma$-orbit) and an element $w\in\W$, and a simple co-weight
$\check{\lambda}$. The elements $\check{\lambda}$ in the dual weight lattice
and $w\in\W$ represent the automorphism $\sigma$ in the following sense: Consider the extended affine Weyl group
$\widehat{W}=\Lambda^{\vee}\rtimes W$ with $\Lambda^{\vee}$ being the co-root lattice of the finite dimensional
Lie algebra $\g$. Then, following \cite{B,T} the stabiliser $ \widehat{W}_A$ of the fundamental alcove $A$,
for its definition see Kac \cite{K}, is naturally isomorphic to $\Lambda^{\vee}/Q^{\vee}$ with
$Q^{\vee}$ being the co-root lattice. Under this identification we have $\sigma=\check{\lambda} w$.
For an explicit description of $\check{\lambda}$ and $w$,
see \cite{T}, Proposition 4.1.4. If $\sigma=id$ we have: 
$w=id$ and $\check{\lambda}=s_r...s_1s_{\theta}(\check{\theta})$,
where $\theta $ is the highest root of the finite dimensional Lie algebra $\g$. 
This leads us to the following decomposition of $\liea $:
$$
\liea=\csa \oplus\bigoplus_{\Orb\subset R}\liea_{\Orb} 
$$
where $\Orb $ runs through all $s_0...s_sw$-orbits of $R$ and $\liea_{\Orb}$ is defined by:
$$
\liea_{\Orb} = \bigoplus_{\alpha\in\Orb}\liea_{\alpha}
$$
Since $\check{\lambda}(z)$ acts trivially on $\csa $ the corresponding $\gamma_{\sigma}$-invariant
subspaces of $\csa $ are direct sums of eigen-spaces of $s_0...s_sw$. These are easily seen to 
provide subbundles of degree zero.
For a given orbit $\Orb $ the loop $\check{\lambda}(z)$ acts diagonally on the root spaces while
$s_0...s_sw$ permutes them cyclically. Therefore, the degree of the subbundle corresponding to 
$\liea_{\Orb}$ is given by $d_{\Orb }=\sum_{\alpha\in\Orb}\alpha(w^{-1}(\check{\lambda}))$.
Denote by $\Gamma $ the group generated by $s_0...s_sw$ and by $\Gamma_{\alpha}$ the 
stabiliser of a given root $\alpha $. Note that the action of $s_0...s_sw$ coincides with the action 
of $\mbox{cox}^{\sigma}$ on $\csa$.
Writing $\check{b}=\sum_{i=1}^{|\Gamma|}(\mbox{cox}^{\sigma})^i(w^{-1}(\check{\lambda}))$
we obtain $d_{\Orb }=\frac{1}{|\Gamma_{\alpha}|}\alpha(\check{b})$ for any $\alpha\in\Orb$.  
The co-weight $\check{b}$ is $\mbox{cox}^{\sigma}$-invariant and non-zero.
Indeed, by \cite{M}, Corollary 2.6 and Proposition 2.7, we know that $\mbox{cox}^{\sigma}$ 
is of infinite order on $\csa\oplus\C C$
the Cartan subalgebra of the centrally extended loop group. Since $s_0...s_sw$ has finite order
the element $(\mbox{cox}^{\sigma})^{|\Gamma|}$ of the extended affine Weyl group has to be a non-vanishing
translation which is easily calculated to be $\check{b}$. The invariance property is clear by definition.
Applying \cite{M} Proposition 2.8 (and calculating $\mbox{mod}\C C$) $\check{b}$ has to coincide,
up to a non-vanishing factor, with the 
solution $b´$ to the equation $(\mbox{cox}^{\sigma}-1)(b´)=kC$.
Thus we obtain a partition of the root system:
\begin{eqnarray}
R_{\check{b}}^{+}&=&\{\alpha\in R|\alpha(\check{b})>0\}\\
R_{\check{b}}^{-}&=&\{\alpha\in R|\alpha(\check{b})<0\}\\
R_{\check{b}}^{0}&=&\{\alpha\in R|\alpha(\check{b})=0\}
\end{eqnarray}
Therefore the root sub-system $R_{\check{b}}^{0}$ describes the root system of some reduction 
Levi subgroup $L(\check{b})=C_{G}(\C \check{b})$ and $R_{\check{b}}^{0}\cup R_{\check{b}}^{+}$
that of the parabolic subgroup $P(\check{b})$. 
If this is not the Harder-Narasimhan reduction, then the Levi subgroup is too big and hence there has to be 
some orbit $\Orb\in R_{\check{b}}^{0}$ and a subspace $\mathfrak{a}\subset\liea_{\Orb }$ such that the bundle
$\mathfrak{A}\subset \mbox{ad}\xi_{\sigma }$
corresponding to $\mathfrak{a}$ has positive degree. This implies that 
$0<\mbox{dim}H^0(E_q,\mathfrak{A})\leq \mbox{dim}H^0(E_q,\mathfrak{G}_{\Orb })$ where $\mathfrak{G}_{\Orb }$
is the bundle corresponding to $\liea_{\Orb }$.
However, in step 2 below we will prove $\mbox{dim}H^0(E_q,\mathfrak{G}_{\Orb })=0$. 

A case-by-case investigation whose basic results are 
summarised in the appendix shows:
\begin{Lemma}
The parabolic subgroup $P(\check{b})$ is maximal and $s_0...s_sw$ is the Coxeter element of its Levi 
$L(\check{b})$. This Levi coincides with the Levi corresponding to an unstable bundle of least possible
automorphism group dimension (see \cite{HS}, Section 6.)
\end{Lemma}

{\bf Step 2:} Next, we have to calculate the dimension of the automorphism group $\mbox{Aut}_G(\xi_{\sigma})$.
According to \cite{HS} Proposition 2.4 this group has a semidirect product structure:
$\mbox{Aut}_G(\xi_{\sigma})=\mbox{Aut}_{L(\check{b})}(\xi_{\sigma L(\check{b})})\ltimes\mbox{Aut}_G(\xi_{\sigma})^{+}$.
The group ${Aut}_G(\xi_{\sigma})^{+}$ is the connected 
unipotent group with Lie algebra 
$\mbox{Lie}\,\mbox{Aut}_G(\xi_{\sigma})^{+}=
H^0(E_q,\xi_{\sigma L(\check{b})}\times^{L(\check{b})}\n)$ where 
$\xi_{\sigma L(\check{b})}\times^{L(\check{b})}\n$ denotes the associated vector bundle
bundle with fibre $\n$ being the Lie algebra of the nilpotent radical of the Harder-Narashiman parabolic subgroup.
The formulae in \cite{HS} preceding Proposition 2.4 of \cite{HS} imply
$$
\mbox{dim}\,\mbox{Aut}_G(\xi_{\sigma})^{+}=\sum_{\Orb\subset R^{+}\backslash R(L(\check{b}))}d_{\Orb}.
$$
These dimensions $\mbox{dim}\,\mbox{Aut}_G(\xi_{\sigma})^{+}$ turn out to be equal to $s$. 
For the readers convenience they are compiled in the appendix.

This leaves us with showing $\mbox{dim}\,\mbox{Aut}_{L(\check{b})}(\xi_{\sigma L(\check{b})})=1$. Using the equality \\
$\mbox{Lie}\,\mbox{Aut}_{L(\check{b})}(\xi_{\sigma L(\check{b})})=H^0(E_q,\mbox{ad}(\xi_{\sigma L(\check{b})}))$ and lifting 
the sections to the trivial $\mbox{Lie}\,L(\check{b})$-bundles over $\Cst$ we have to determine
the solutions $h\in\L(\mbox{Lie}\,L(\check{b}))$ of the equation
$$
\mbox{Ad}(\gamma_{\sigma}(z))(h(z))=h(qz).
$$
Without loss of generality we can assume that $h$ takes values either in an eigen-space of 
the Coxeter element $\mbox{cox}_L$ or in a space $\liea_{\Orb }$ with $\Orb\subset R(L(\check{b}))$.
In the first case we are left with the equation $\zeta h(z)=h(qz)$ where $\zeta $ is the 
corresponding eigenvalue. The only non-vanishing solution of this equation occur for
$\zeta=1$ and constant functions on $E_q$. The corresponding eigen-space coincides with the 
one-dimensional centre of $\mbox{Lie}\,L(\check{b})$.
Labelling the elements of the orbit $\Orb=\{\beta_1,...,\beta_p\}$ for $p=|\Orb |$ and 
setting $h(z)=\sum_{i=1}^{p}h_i(z)\otimes x_{\beta_i}$ where $x_{\beta_i}$ is a basis element 
of the corresponding root space we require:
$$
z^{(\beta_i,w^{-1}(\check{\lambda}))}h_i(z)=h_{i+1}(qz).
$$
Iterating this formula we obtain:

$$
h_1(q^pz)=q^{\sum_{i=1}^{p}(p-i)(\beta_i,w^{-1}(\check{\lambda}))}z^{d_{\Orb}}h_{1}(z).
$$
A tedious calculation using $d_{\Orb}=0$ shows that the exponent of $q$ on the right hand side is 
not divisible by $p$. Therefore, this functional equation does not permit a solution on $\Cst $. 
\end{proof}
\begin{Remark} 
(i) This case-by-case line of reasoning works for simply connected structure groups and similar 
phenomena occur (the role of $\check{b}$, the restriction of the Coxeter element of the affine group yielding 
the Coxeter element of the Levi from the reduction, etc.). We refrain from giving the detail in favour of the 
more uniform and elegant treatment we use in this case.\\
(ii)
It would be interesting to find a more
conceptual reason of the fact that the restriction of the twisted Coxeter
element to $\csa $ gives the Coxeter element of the Levi $L(\check{b})$.
\end{Remark}



\section{Appendix}

Here we summarise the explicit results of the case-by-case calculations.
Let us label the vertexes of the Dynkin diagram as in \cite{B}. 
For the diagram automorphisms we use the following conventions.
In all cases except $\mbox{D}_{2n}$, $\gamma $ is the generator of
the $\pi_1(G^{ad})$, where $G^{ad}$ is the adjoint group. 
In the $D_{2n}$-case, $\gamma^2$ generates the fundamental group of 
$SO_{4n}$ and $\tau\neq id$ generates $\pi_1(G)$ with $\tau^2=id$ and 
$G$ is not isomorphic to $SO_{4n}$.
We only indicate $d_{\Orb}$ for orbits containing the basis elements because $d_{\Orb}$ for all other
orbits can be derived from these. (Note, that taking the sum over elements of two orbits might change the
orbit length!) In the sequel $\alpha_0$ denotes the negative highest root of the root system 
$\Delta$ of $G$ corresponding to the basis $\Pi=\{\alpha_1,\ldots,\alpha_r\}$.

{\bf $\mbox{A}_n:$} For $\sigma=id $ the $\mbox{cox}$-orbits involving the basis elements
look as follows:\\
$\Orb_1=\{\alpha_1,...,\alpha_{n-1},-(\alpha_1+...+\alpha_{n-1})\}$ and
$\Orb_2=\{\alpha_n, \alpha_{n-1}+\alpha_n,...,\alpha_1+...+\alpha_n\}$.\\
Evaluation on $\check{b}$ yields: $d_{\Orb_1}=0$ and $d_{\Orb_2}=n+1$.
Hence, $P(\check{b})=P_{\alpha_n}$ and $\mbox{dim}\,\mbox{Aut}_G(\xi_{\sigma})^{+}=n+1$

For $\sigma=\gamma $ we have $\check{\lambda}=\check{\lambda}_1$ and $w=s_1...s_n$.
Thus, $\mbox{cox}^{\sigma}=s_1\check{\lambda}_1s_1...s_n$. 
Considering the $\mbox{cox}^{\sigma}$-orbits involving the basis elements yields:\\
$\Orb_1=\{\alpha_2,...,\alpha_n,-(\alpha_2+...+\alpha_n)\}$ and
$\Orb_2=\{\alpha_1, \alpha_1+\alpha_2,...,\alpha_1+...+\alpha_n\}$.\\ 
Evaluation on $\check{b}$ yields: $d_{\Orb_1}=0$ and $d_{\Orb_2}=-1$.
Hence, $P(\check{b})=P_{-\alpha_1}$ and $\mbox{dim}\,\mbox{Aut}_G(\xi_{\sigma})^{+}=1$

For $\sigma=\gamma^l $ with $l|n$ we have $\mbox{cox}^{\sigma}=s_1...s_l\check{\lambda}_lw^l$.
We obtain the same results as above except for $d_{\Orb_2}$ which yields $d_{\Orb_2}=-l$
and $\mbox{dim}\,\mbox{Aut}_G(\xi_{\sigma})^{+}=l$.

{\bf $\mbox{B}_n:$} Consider $\sigma=id $. 
For displaying the correct parabolic subgroup in standard form it turns out that 
a change of the basis of the root system is necessary:\\
$\beta_i=\alpha_{n-i}$, $1\leq i\leq n-1$ and
$\beta_n=-(\alpha_1+...+\alpha_n)$.\\
Then, the $\mbox{cox}$-orbits involving the basis elements have the following form:\\
$\Orb_1=\{\beta_1,...,\beta_{n-2},-(\beta_1+...+\beta_{n-2})\}$,
$\Orb_2=\{\beta_n,-\beta_n\}$ and for $n$ odd
$\Orb_3=\{\beta_{2i}+...+\beta_{n-1}, 1\leq i\leq \frac{n-1}{2}\}
\cup\{\beta_{2i-1}+...+\beta_{n-1}+2\beta_n, 1\leq i\leq \frac{n-1}{2}\}$,
respectively for $n$ even
$\Orb_3=\{\beta_{i}+...+\beta_{n-1}, 1\leq i\leq n-1\}
\cup\{\beta_{i}+...+\beta_{n-1}+2\beta_n, 1\leq i\leq n-1\}$.\\
Evaluation on $\check{b}$ yields: $d_{\Orb_1}=0$, $d_{\Orb_2}=0$ and $d_{\Orb_3}=-1$,
for odd $n$ respectively $d_{\Orb_3}=-2$ in the even case.
Hence, $P(\check{b})=P_{-\beta_{n-1}}$ and $\mbox{dim}\,\mbox{Aut}_G(\xi_{\sigma})^{+}=n+1$.

Consider $\sigma=\gamma $. We have $\check{\lambda}=\check{\lambda}_1$ while $w$
fixes $\alpha_2,...,\alpha_n$ and interchanges $\alpha_1$ with $\alpha_0$.
Hence, $\mbox{cox}^{\sigma}=s_n...s_1\check{\lambda}_1w$.
The $\mbox{cox}^{\sigma}$-orbits involving the basis elements have the following shape:\\
$\Orb_1=\{\alpha_1,...,\alpha_{n-1},-(\alpha_1+...+\alpha_{n-1})\}$ and
$\Orb_2=\{\alpha_n, \alpha_{n-1}+\alpha_n,...,\alpha_1+...+\alpha_n\}$.\\ 
Evaluation on $\check{b}$ yields: $d_{\Orb_1}=0$ and $d_{\Orb_2}=-1$.
Hence, $P(\check{b})=P_{-\alpha_n}$ and $\mbox{dim}\,\mbox{Aut}_G(\xi_{\sigma})^{+}=n$.

{\bf $\mbox{C}_n:$} For $\sigma=id $ the $\mbox{cox}$-orbits involving the basis 
elements look as follows:\\ 
$\Orb_1=\{\alpha_1,...,\alpha_{n-1},-(\alpha_1+...+\alpha_{n-1})\}$,
$\Orb_2=\{2(\alpha_1+...+\alpha_{n-1})+\alpha_n,....,2\alpha_{n-1}+\alpha_n, \alpha_n\}$.\\
Evaluation on $\check{b}$ yields: $d_{\Orb_1}=0$, $d_{\Orb_2}=2$,
Hence, $P(\check{b})=P_{\alpha_{n}}$ and $\mbox{dim}\,\mbox{Aut}_G(\xi_{\sigma})^{+}=n+1$.

For $\sigma=\gamma $ we have $\check{\lambda}=\check{\lambda}_n$ and $w$ interchanging
$\alpha_i$ with $\alpha_{n-i}$ with $\alpha_0$ being the negative of the highest root of the finite dimensional Lie algebra $\g $.
Thus, $\mbox{cox}^{\sigma}=s_{[\frac{n+1}{2}]}...s_n1\check{\lambda}_nw$. 
For displaying the correct parabolic subgroup in standard form it turns out that 
a change of the basis of the root system is necessary:\\
$\beta_i=\alpha_{n-i}$, $1\leq i\leq [\frac{n-1}{2}]-1$, 
$\beta_{[\frac{n-1}{2}]}=-(\alpha_1+...\alpha_{[\frac{n}{2}]+1}+2(\alpha_{[\frac{n}{2}]+2}+...+\alpha_{n-1})+\alpha_n)$
$\beta_i=\alpha_{i-[\frac{n-1}{2}]}$ for $[\frac{n+1}{2}]\leq i\leq n-2$,
$\beta_n=\alpha_{[\frac{n}{2}]}$ and
$\beta_n=2(\alpha_{[\frac{n}{2}]+1}+...+\alpha_{n-1})+\alpha_n$
 
For even $n$ the $\mbox{cox}^{\sigma}$-orbits involving the basis elements are given by:\\
$\Orb_1=\{\beta_1,...,\beta_{n-2},-(\beta_1+...+\beta_{n-2})\}$,
$\Orb_2=\{\beta_{n},-\beta_{n}\}$ and 
$\Orb_3=\{\beta_i+...+\beta_{n-1}, 1\leq i\leq n-1\}\cup\{\beta_i+...+\beta_{n}, 1\leq i\leq n-1\}$.\\
Evaluation on $\check{b}$ yields: $d_{\Orb_1}=0$, $d_{\Orb_2}=0$ and $d_{\Orb_3}=-1$.
Hence, $P(\check{b})=P_{-\beta_{n-1}}$ and $\mbox{dim}\,\mbox{Aut}_G(\xi_{\sigma})^{+}=[\frac{n}{2}]-1$.

In the odd case we calculate:\\
$\Orb_1=\{\beta_1,...,\beta_{n-1},-(\beta_1+...+\beta_{n-1})\}$ and
$\Orb_2=\{2(\beta_i,...,\beta_{n-1})+\beta_n, 1\leq i\leq n\}$.\\
The degrees are given by $d_{\Orb_1}=0$ and $d_{\Orb_2}=-1$ yielding $P(\check{b})=P_{-\beta_{n}}$.

{\bf $\mbox{D}_n:$} Consider $\sigma=id $. Also, here we have to introduce a new basis of the
root system:\\
$\beta_i=\alpha_{n-i-1}$, $1\leq i\leq n-2$, 
$\beta_{n-1}=-(\alpha_1+...+\alpha_{n-1})$and
$\beta_n=-(\alpha_1+...+\alpha_{n-2}+\alpha_n)$.\\
Then, the $\mbox{cox}$-orbits involving the basis elements have the following form:\\
$\Orb_1=\{\beta_1,...,\beta_{n-3},-(\beta_1+...+\beta_{n-3})\}$,
$\Orb_2=\{\beta_{n-1},-\beta_{n-1}\}$,
$\Orb_3=\{\beta_n,-\beta_n\}$ and for even $n$ 
$\Orb_4=\{\beta_{2i}+...+\beta_{n-2}, 1\leq i\leq \frac{n-2}{2}\}
\cup\{\beta_{2i-1}+...+\beta_{n-1}+\beta_n, 1\leq i\leq \frac{n-2}{2}\}$,
respectively for $n$ odd
$\Orb_4=\{\beta_{i}+...+\beta_{n-2}, 1\leq i\leq n-1\}
\cup\{\beta_{i}+...+\beta_{n-1}+\beta_n, 1\leq i\leq n-1\}$.\\
Evaluation on $\check{b}$ yields: $d_{\Orb_1}=d_{\Orb_2}=d_{\Orb_3}=0$ and $d_{\Orb_4}=-1$,
for even $n$ respectively $d_{\Orb_4}=-2$ in the odd case.
Hence, $P(\check{b})=P_{-\beta_{n-2}}$ and $\mbox{dim}\,\mbox{Aut}_G(\xi_{\sigma})^{+}=n+1$.

For $\sigma=\gamma^2 $ i.e. $G=SO_{2n}$ we have $\check{\lambda}=\check{\lambda}_1$ and $w$
simultaneously interchanges $\alpha_0$ with $\alpha_1$ and $\alpha_{n-1}$ with $\alpha_n$.
Thus, $\mbox{cox}^{\sigma}=s_{n-1}...s_1\check{\lambda}_1w$. The $\mbox{cox}^{\sigma}$-orbits involving the
basis elements have the following shape:\\
$\Orb_1=\{\alpha_1,...,\alpha_{n-2},\alpha_n ,-(\alpha_1+...+\alpha_{n-2}+\alpha_n)\}$ and
$\Orb_2=\{\alpha_i+2(\alpha_{i+1}+...+\alpha_{n-2})+\alpha_{n-1}+\alpha_{n}, 1\leq i\leq n-2\}
\cup\{\alpha_{n-1},\alpha_1+...+\alpha_{n-1}\}$.\\ 
Evaluation on $\check{b}$ yields: $d_{\Orb_1}=0$ and $d_{\Orb_2}=-2$.
Hence, $P(\check{b})=P_{-\alpha_{n-1}}$ and $\mbox{dim}\,\mbox{Aut}_G(\xi_{\sigma})^{+}=n-1$.

Consider $n$ odd and $\sigma=\gamma$. Then $\check{\lambda}=\check{\lambda}_n$ and $w$ permutes
$\alpha_0,\alpha_n,\alpha_1$ and $\alpha_{n-1}$ cyclically while interchanging $\alpha_i$
with $\alpha_{n-i}$ for the other labels. We have 
$\mbox{cox}^{\sigma}=s_{\frac{n+1}{2}}...s_n\check{\lambda}_nw$. Again we have to find a new basis in 
order to get a parabolic in standard form:\\
$\beta_i=\alpha_{i+1}$ for $1\leq i\leq \frac{n-1}{2}-1$,
$\beta_{\frac{n-1}{2}}=\alpha_{\frac{n+1}{2}}+...+\alpha_n$,
$\beta_i=\alpha_{\frac{3n-1}{2}-i}$, $\frac{n+1}{2}\leq i\leq n-2$,
$\beta_{n-1}=\alpha_1+...+\alpha_{\frac{n+1}{2}}$ and
$\beta_n=-(\alpha_1+...+\alpha_{\frac{n+1}{2}}+2(\alpha_{\frac{n+1}{2}+1}+...+\alpha_{n-2})+\alpha_{n-1}+\alpha_n)$.\\
The orbits containing simple roots look like:\\
$\Orb_1=\{\beta_1,...,\beta_{n-1},-(\beta_1+...+\beta_{n-1})\}$
$\Orb_2=\{\beta_i+2(\beta_{i+1}+...+\beta_{n-2})+\beta_{n-1}+\beta_{n}, 1\leq i\leq n-2\}
\cup\{\beta_n,\beta_1+...+\beta_{n-2}+\beta_{n}\}$\\
For the degrees we calculate $d_{\Orb_1}=0$ and $d_{\Orb_2}=2$ yielding $P(\check{b})=P_{-\beta_{n}}$
and $\mbox{dim}\,\mbox{Aut}_G(\xi_{\sigma})^{+}=\frac{n-1}{2}$.

Let us turn to even $n$ and $\sigma=\tau $. Then, $\check{\lambda}=\check{\lambda}_n$ and
$w$ interchanges $\alpha_i $ and $\alpha_{n-i}$ while $\mbox{cox}^{\sigma}=s_{\frac{n}{2}}...s_n\check{\lambda}_nw$.
Again we have to find a new basis in order to get a parabolic in standard form:\\
$\beta_i=\alpha_{i+1+\frac{n}{2}}$ for $1\leq i\leq \frac{n}{2}-3$,
$\beta_{\frac{n}{2}-2}=\alpha_{\frac{n}{2}}+...+\alpha_n$,
$\beta_i=\alpha_{n-i-2}$ for $\frac{n}{2}-1\leq i\leq n-4$,
$\beta_{n-3}=\alpha_1$, 
$\beta_{n-2}=-(\alpha_1+...+\alpha_{\frac{n}{2}})$, 
$\beta_{n-1}=-(\alpha_{\frac{n}{2}+1}+...+\alpha_{n-1})$ and 
$\beta_{n}=-(\alpha_{\frac{n}{2}+1}+...+\alpha_{n-2}+\alpha_{n})$.\\
We obtain the following orbits containing simple roots:\\
$\Orb_1=\{\beta_1,...,\beta_{n-4},-(\beta_1+...+\beta_{n-4})\}$
$\Orb_2=\{\beta_{n-2},\beta_{n-1},\beta_{n},-(\beta_{n-2}+\beta_{n-1}+\beta_{n})\}$ and
$\Orb_3=\{\beta_i+...\beta_{n-3}+r_j, 1\leq i\leq n-3, r_j
\in\{0,\beta_{n-2}+\beta_{n-1},\beta_{n-2}+\beta_n,\beta_{n-2}+\beta_{n-1}+\beta_n,\}\}$.\\
The degrees are given by $d_{\Orb_1}=d_{\Orb_2}=0$ and 
$d_{\Orb_3}=2$ yielding $P(\check{b})=P_{\beta_{n-3}}$
and $\mbox{dim}\,\mbox{Aut}_G(\xi_{\sigma})^{+}=\frac{n}{2}+1$.

{\bf $\mbox{E}_6$:} For $\sigma=id $ we have to find a new basis of the root system:\\
$\beta_1=-(\alpha_1+2\alpha_2+2\alpha_3+2\alpha_4+\alpha_5+\alpha_6)$, 
$\beta_2=\alpha_2+\alpha_3+\alpha_4+\alpha_5$,
$\beta_3=\alpha_1$, $\beta_4=\alpha_2+\alpha_3+\alpha_6$, 
$\beta_5=\alpha_3+\alpha_4$ and
$\beta_6=-(\alpha_1+\alpha_2+2\alpha_3+\alpha_4+\alpha_5+\alpha_6)$.\\
Then, the $\mbox{cox}$-orbits involving the basis elements have the following form:\\
$\Orb_1=\{\beta_1,\beta_{2},-\beta_1-\beta_{2}\}$,
$\Orb_2=\{\beta_4,\beta_{5},-\beta_4-\beta_{5}\}$,
$\Orb_3=\{\beta_6,-\beta_6\}$ and
$\Orb_4=\{\beta_3,\beta_3+\beta_6,\beta_2+\beta_3+\beta_4,\beta_2+\beta_3+\beta_4+\beta_6,
 \beta_1+\beta_2+\beta_3+\beta_4+\beta_5,\beta_1+\beta_2+\beta_3+\beta_4+\beta_5+\beta_6\}$.\\
Evaluation on $\check{b}$ yields: $d_{\Orb_1}=d_{\Orb_2}=d_{\Orb_3}=0$ and $d_{\Orb_4}=-1$.
Hence, $P(\check{b})=P_{-\beta_{3}}$ and $\mbox{dim}\,\mbox{Aut}_G(\xi_{\sigma})^{+}=7$.

Consider $\sigma=\gamma$. Then, $\check{\lambda}=\check{\lambda}_6$ and $w$
permutes $\alpha_0,\alpha_1,\alpha_6$ respectively $\alpha_2,\alpha_3,\alpha_5$ cyclically
while fixing $\alpha_4$. The twisted Coxeter element has the form 
$s_1s_3s_4\check{\lambda}_6w$.\\
Also here we need  new basis of the root system:\\
$\beta_1=\alpha_5+\alpha6$, $\beta_2=\alpha_1+\alpha_3+\alpha_4+\alpha_5$,
$\beta_3=\alpha_3+\alpha4$, $\beta_4=\alpha_2$, 
$\beta_5=\alpha_0=-(\alpha_1+2\alpha_2+2\alpha_3+3\alpha_4+2\alpha_5+\alpha_6)$ and
$\beta_6=\alpha_1+\alpha_2+\alpha_3+2\alpha_4+\alpha_5+\alpha_6$.\\
The orbits containing simple roots look as follows:\\
$\Orb_1=\{\beta_1,\beta_2,\beta_3,\beta_4,-(\beta_1+\beta_2+\beta_3+\beta_4)\}$,
$\Orb_2=\{\beta_6,-\beta_6\}$ and
$\Orb_3=\{\beta_5,\beta_2+\beta_4+\beta_5+\beta_6,\beta_1+\beta_3+\beta_4+\beta_5+\beta_6, 
\beta_1+\beta_2+2\beta_3+2\beta_4+\beta_5, \beta_2+\beta_3+2\beta_4+\beta_5+\beta_6,
\beta_3+\beta_4+\beta_5,\beta_5+\beta_6,\beta_1+\beta_3+\beta_4+\beta_5, 
\beta_1+\beta_2+2\beta_3+2\beta_4+\beta_5+\beta_6,\beta_2+\beta_3+2\beta_4+\beta_5\}$.\\
The degrees are given by $d_{\Orb_1}=d_{\Orb_2}=0$ and $d_{\Orb_3}=-1$. 
Hence $P(\check{b})=P_{-\beta_{5}}$ and $\mbox{dim}\,\mbox{Aut}_G(\xi_{\sigma})^{+}=3$.

{\bf $\mbox{E}_7$:} Consider $\sigma=id $. We have to find a new basis of the root system:\\
$\beta_1=-(\alpha_1+2\alpha_2+2\alpha_3+2\alpha_4+\alpha_5+\alpha_7)$, 
$\beta_2=\alpha_2+\alpha_3+\alpha_4+\alpha_5$,
$\beta_3=\alpha_1$, $\beta_4=\alpha_2+\alpha_3+\alpha_7$, 
$\beta_5=\alpha_3+\alpha_4$,
$\beta_6=-(\alpha_1+2\alpha_2+3\alpha_3+2\alpha_4+\alpha_5+\alpha_6+\alpha_7)$ and
$\beta_6=-(\alpha_1+\alpha_2+2\alpha_3+\alpha_4+\alpha_5+\alpha_7)$.\\
Then, the $\mbox{cox}$-orbits involving the basis elements have the following form:\\
$\Orb_1=\{\beta_1,\beta_{2},-\beta_1-\beta_{2}\}$,
$\Orb_2=\{\beta_4,\beta_{5},\beta_{6},-\beta_4-\beta_{5}-\beta_{6}\}$,
$\Orb_3=\{\beta_7,-\beta_7\}$ and
$\Orb_4=\{\beta_3, \beta_2+\beta_3,\beta_2+\beta_3+\beta_4+\beta_7,\beta_1+\beta_2+\beta_3+\beta_4+\beta_5,
\beta_3+\beta_4+\beta_5+\beta_6+\beta_7,\beta_1+\beta_2+\beta_3+\beta_4+\beta_7,\beta_3+\beta_4+\beta_5,
\beta_2+\beta_3+\beta_4+\beta_5+\beta_6+\beta_7,\beta_1+\beta_2+\beta_3,\beta_3+\beta_4+\beta_7, 
\beta_2+\beta_3+\beta_4+\beta_5,\beta_1+\beta_2+\beta_3+\beta_4+\beta_5+\beta_6+\beta_7\}$.\\
Evaluation on $\check{b}$ yields: $d_{\Orb_1}=d_{\Orb_2}=d_{\Orb_3}=0$ and $d_{\Orb_4}=-1$.
Hence, $P(\check{b})=P_{-\beta_{3}}$ and $\mbox{dim}\,\mbox{Aut}_G(\xi_{\sigma})^{+}=8$.

For $\sigma=\gamma $ we get $\check{\lambda}=\check{\lambda}_7$ and $w$
interchanging $\alpha_0$ with $\alpha_7$, $\alpha_1$ with $\alpha_6$ and $\alpha_3$ with $\alpha_5$
fixing $\alpha_2$ and $\alpha_4$. Thus $\mbox{cox}^{\sigma}=s_7s_6s_5s_4s_2\check{\lambda}_7w$.
A basis in which the parabolic subgroup will be of standard type looks as follows:\\
$\beta_1=\alpha_2+\alpha_4+\alpha_5+\alpha_6$,
$\beta_2=\alpha_2+\alpha_3+\alpha_4+\alpha_5+\alpha_6+\alpha_7$,
$\beta_3=\alpha_1+\alpha_3$, $\beta_4=\alpha_4+\alpha_5$, $\beta_5=-\alpha_5$,
$\beta_6=-(\alpha_1+\alpha_2+\alpha_3+2\alpha_4+\alpha_5+\alpha_6+\alpha_7)$ and
$\beta_7=-(\alpha_3+\alpha_4+\alpha_5+\alpha_6)$.\\
The orbits containing simple roots are:\\
$\Orb_1=\{\beta_1,\beta_2,\beta_3,\beta_{4},-(\beta_1+\beta_2+\beta_3+\beta_{4})\}$,
$\Orb_2=\{\beta_6,\beta_7,-(\beta_6+\beta_7)\}$ and 
$\Orb_3=\{\beta_5,\beta_1+\beta_3+\beta_4+\beta_5+\beta_6+\beta_7,\beta_1+\beta_2+2\beta_3+2\beta_4+\beta_5+\beta_6, 
\beta_2+\beta_3+2\beta_4+\beta_5,\beta_2+\beta_4+\beta_5+\beta_6+\beta_7, \beta_5+\beta_6, 
\beta_1+\beta_3+\beta_4+\beta_5,\beta_1+\beta_2+2\beta_3+2\beta_4+\beta_5+\beta_6+\beta_7\}$.\\
We calculate for the degrees $d_{\Orb_1}=d_{\Orb_2}=0$ and $d_{\Orb_3}=1$ implying $P(\check{b})=P_{\beta_{5}}$
and $\mbox{dim}\,\mbox{Aut}_G(\xi_{\sigma})^{+}=5$.

{\bf $\mbox{E}_8$:} For $\sigma=id $ a new basis of the root system is given by:\\
$\beta_1=-(\alpha_1+2\alpha_2+3\alpha_3+4\alpha_4+4\alpha_5+3\alpha_6+\alpha_7+2\alpha_8)$, 
$\beta_2=\alpha_4+\alpha_5+\alpha_6$,
$\beta_3=\alpha_3+\alpha_4+\alpha_5+\alpha_8$,
$\beta_4=\alpha_2+\alpha_3+\alpha_4+\alpha_5+\alpha_6+\alpha_7$
$\beta_5=\alpha_1$,  
$\beta_6=\alpha_2+\alpha_3+\alpha_4+2\alpha_5+\alpha_6+\alpha_8)$,
$\beta_7=-(\alpha_1+2\alpha_2+2\alpha_3+3\alpha_4+4\alpha_5+2\alpha_6+\alpha_7+2\alpha_8)$ and
$\beta_7=-(\alpha_1+\alpha_2+2\alpha_3+2\alpha_4+3\alpha_5+2\alpha_6+\alpha_7+\alpha_8)$.\\
Then, the $\mbox{cox}$-orbits involving the basis elements have the following form:\\
$\Orb_1=\{\beta_1,\beta_{2},\beta_3,\beta_4,-\beta_1-\beta_{2}-\beta_3-\beta_4\}$,
$\Orb_2=\{\beta_6,\beta_{7},-\beta_{6}-\beta_{7}\}$,
$\Orb_3=\{\beta_8,-\beta_8\}$ and
$\Orb_4=\{\beta_i...\beta_5, \beta_i...\beta_5+\beta_6,\beta_i...\beta_5+\beta_6+\beta_7,
\beta_i...\beta_5+\beta_8,\beta_i...\beta_5+\beta_6+\beta_8,
\beta_i...\beta_5+\beta_6+\beta_7+\beta_8, 1\leq i\leq 5\}$.\\
Evaluation on $\check{b}$ yields: $d_{\Orb_1}=d_{\Orb_2}=d_{\Orb_3}=0$ and $d_{\Orb_4}=-1$.
Hence, $P(\check{b})=P_{-\beta_{5}}$ and $\mbox{dim}\,\mbox{Aut}_G(\xi_{\sigma})^{+}=9$.  

{\bf $\mbox{F}_4$:} For $\sigma=id $ a new basis of the root system has to be introduced:\\
$\beta_1=-(\alpha_1+\alpha_2+2\alpha_3)$, 
$\beta_2=\alpha_1$,
$\beta_3=\alpha_2+\alpha_3$, 
$\beta_4=-(\alpha_1+2\alpha_2+2\alpha_3+\alpha_4)$.\\
Then, the $\mbox{cox}$-orbits involving the basis elements have the following form:\\
$\Orb_1=\{\beta_1,-\beta_1\}$,
$\Orb_2=\{\beta_3,\beta_{4},-\beta_{3}-\beta_{4}\}$,
$\Orb_3=\{\beta_2, \beta_1+\beta_2, \beta_2+2\beta_3, \beta_2+\beta_3+2\beta_4,\beta_1+\beta_2+2\beta_3 
\beta_1+\beta_2+2\beta_3+2\beta_4\}$.\\
Evaluation on $\check{b}$ yields: $d_{\Orb_1}=d_{\Orb_2}=0$ and $d_{\Orb_3}=-1$.
Hence, $P(\check{b})=P_{-\beta_{2}}$ and $\mbox{dim}\,\mbox{Aut}_G(\xi_{\sigma})^{+}=5$.  

{\bf $\mbox{G}_2$:} Consider $\sigma=id $. Introduce a new basis of the root system:\\
$\beta_1=\alpha_1$, 
$\beta_2=-\alpha_1-\alpha_2$.\\
Then, the $\mbox{cox}$-orbits involving the basis elements have the following form:\\
$\Orb_1=\{\beta_2,-\beta_2\}$,
$\Orb_3=\{\beta_1, \beta_1+3\beta_2\}$.\\
Evaluation on $\check{b}$ yields: $d_{\Orb_1}=0$ and $d_{\Orb_3}=-1$.
Hence, $P(\check{b})=P_{-\beta_{1}}$ and $\mbox{dim}\,\mbox{Aut}_G(\xi_{\sigma})^{+}=3$.



\end{document}